\documentclass[11pt]{Mytran-l}
\usepackage[ansinew]{inputenc}
\usepackage[centertags]{amsmath}
\usepackage{exscale}
\usepackage{relsize}
\usepackage{amsfonts}
\usepackage{amssymb}
\usepackage{amsthm}
\usepackage{newlfont}
\usepackage{color}
\usepackage[colorlinks,citecolor=blue,urlcolor=blue,breaklinks=true]{hyperref}

\mathchardef\ordinarycolon\mathcode`\:
  \mathcode`\:=\string"8000
  \begingroup \catcode`\:=\active
    \gdef:{\mathrel{\mathop\ordinarycolon}}
  \endgroup

\newtheorem*{mthm*}{Main Theorem}

\newtheorem{thm}{Theorem}
\newtheorem{lem}[thm]{Lemma}
\newtheorem{prop}[thm]{Proposition}

\newtheorem{examp}[thm]{Example}

\newcommand{\D}{\mathcal{D}}
\newcommand{\B}{\mathcal{B}}

\newcommand{\h}{\mathcal{H}}
\newcommand{\aut}{\emph{Aut}}
\newcommand{\Aut}{\mbox{Aut}}

\begin{document}

\title[Almost simple groups with socle $L_n(q)$ acting on SQS$(v)$]
{Almost simple groups with socle $L_n(q)$ acting on Steiner quadruple systems}

\author{Michael Huber}

\address{Wilhelm Schickard Institute for Computer Science, University of Tuebingen, Sand~13,
D-72076 Tuebingen, Germany}

\email{michael.huber@uni-tuebingen.de}



\thanks{The author gratefully acknowledges support by the Deutsche
Forschungsgemeinschaft (DFG) via a Heisenberg grant (Hu954/4).}

\date{January 16, 2009; and in revised form July 7, 2009}

\commby{The Managing Editor}

\maketitle

\begin{abstract}
Let $N=L_n(q)$, \mbox{$n \geq 2$}, $q$ a prime power, be a projective linear simple group.
We classify all Steiner quadruple systems admitting a group $G$ with $N \leq G \leq \Aut(N)$.
In particular, we show that $G$ cannot act as a group of automorphisms on any Steiner quadruple system for $n>2$.
\end{abstract}


\section{Introduction}\label{intro}

All Steiner quadruple systems which admit a flag-transitive group of automorphisms were classified in~\cite{Hu2001}.
The most interesting examples which occur have an almost simple group with socle $L_2(q)$ as group of automorphisms.
In this note, we examine for these type of groups the general case when flag-transitivity is omitted.
Before stating the main result, we present the examples that arise in our consideration.

\begin{examp}\label{Ex1}\em
$\mathcal{D}$ is a Steiner quadruple system \mbox{SQS$(3^d +1)$}
whose points are the elements of the projective line \mbox{$\mathbb{F}_{3^d} \cup \{\infty\}$}
and whose blocks are the images of $\mathbb{F}_3 \cup \{\infty\}$ under
$PGL_2(3^d)$ with $d \geq 2$ (resp. $L_2(3^d)$ with $d
>1$ odd),  and $L_2(3^d) \leq G \leq P \mathit{\Gamma} L_2(3^d)$. The derived design at any given
point is the Steiner triple system \mbox{STS$(3^d)$} whose
points and blocks are the points and lines of the affine space $AG(d,3)$.
In this case, $G$ acts flag-transitively on $\mathcal{D}$.
\end{examp}

\begin{examp}\label{Ex2}\em
$\mathcal{D}$ is a Steiner quadruple system \mbox{SQS$(q+1)$} whose
points are the elements of  \mbox{$\mathbb{F}_q \cup \{\infty\}$} with a
prime power $q \equiv 7$ $(\mbox{mod}$ $12)$ and whose blocks are the
images of $\{0,1,\infty,\varepsilon\}$ under $L_2(q)$, where
$\varepsilon$ is a primitive sixth root of unity in $\mathbb{F}_q$, and
$L_2(q) \leq G \leq P \mathit{\Sigma} L_2(q)$.
The derived design at any given point is the \emph{Netto triple system} $N(q)$,
a detailed description of which can be found in~\cite[Section\,3]{Deletal86}.
Here, $G$ acts flag-transitively on $\mathcal{D}$.
\end{examp}

\begin{examp}\label{Ex3}\em
$\D$ is a Steiner quadruple system \mbox{SQS$(3^{2d}+1)$} whose
points are the elements of \mbox{$\mathbb{F}_{3^{2d}} \cup \{\infty\}$} and
whose blocks are the disjoint union of the images of
$\{0,1,-1,\infty\}$ and $\{0,1,a,\infty\}$
under $L_2(3^{2d})$ with $d \geq 1$, $a \not \in
(\mathbb{F}_{3^{2d}}^*){^2}$, and $L_2(3^{2d}) \leq G \leq P \mathit{\Sigma} L_
2(3^{2d})$.
In this case, $G$ has two orbits on the \mbox{$3$-subsets} and hence on the blocks of $\mathcal{D}$. Therefore, flag-transitivity cannot hold.
\end{examp}

Our main result is as follows.

\begin{mthm*}\label{PSL_SQS}
Let $\mathcal{D}$ be a non-trivial Steiner quadruple system \mbox{SQS$(v)$} of order $v$, and $N \leq G \leq \aut(N)$ with a projective linear simple group $N=L_n(q)$, $n \geq 2$, $q$ a prime power, $v=\frac{q^n-1}{q-1}$.
Then $G \leq \aut(\mathcal{D})$ acts on $\mathcal{D}$ if and only if one of the cases described in Examples~\emph{\ref{Ex1},~\ref{Ex2},~\ref{Ex3}} above occurs (up to isomorphism).
\end{mthm*}


\section{Preliminaries}\label{Prelim}

For positive integers $t \leq k \leq v$ and $\lambda$, we define a
\mbox{\emph{$t$-$(v,k,\lambda)$ design}} to be a finite incidence
structure \mbox{$\mathcal{D}=(X,\mathcal{B},I)$}, where $X$ denotes
a set of \emph{points}, $\left| X \right| =v$, and $\mathcal{B}$ a
set of \emph{blocks}, $\left| \mathcal{B} \right| =b$, satisfying
the following properties: (i) each block $B \in \mathcal{B}$ is
incident with $k$ points; and (ii) each \mbox{$t$-subset} of $X$ is
incident with $\lambda$ blocks. A \emph{flag} is an incident point-block pair.

For historical reasons, a \mbox{$t$-$(v,k,1)$ design} is called a \emph{Steiner \mbox{$t$-design}} or a \emph{Steiner system}.
We note that in this case each block is determined by the set of points which are
incident with it, and thus can be identified with a \mbox{$k$-subset} of $X$ in a unique way.
A \emph{Steiner triple system} of order $v$ is a \mbox{$2$-$(v,3,1)$ design}.
A \emph{Steiner quadruple system} of order $v$ is a \mbox{$3$-$(v,4,1)$ design},
and will be denoted in the following by \mbox{SQS$(v)$}. The case $v=4$ yields the \emph{trivial} \mbox{SQS$(v)$}.
A simple example of a Steiner quadruple  system is the \mbox{SQS$(2^n)$} consisting of the
points and planes of the \mbox{$n$-dimensional} binary affine space $AG(n,2)$ for
each $n\geq 2$. Using recursive constructions, H.~Hanani~\cite{Han1960} showed that
the following condition for the existence of a \mbox{SQS$(v)$} (the necessity of which is easy to see)
is also sufficient:

\begin{prop}{\em (Hanani,~1960).}\label{Hanani}
A Steiner quadruple system SQS$(v)$ of order $v$ exist if
and only if
\[v \equiv 2 \; \mbox{or} \; 4 \; (\emph{mod} \; 6)\;\, (v \geq 4).\]
\end{prop}

For $v=8$ and $v=10$ there exists a \mbox{SQS$(v)$} in each case, unique up to isomorphism. These are
the affine space $AG(3,2)$ and the M\"{o}bius plane of order $3$, see Barrau~\cite{Bar1908}, 1908. For $v=14$ we have exactly four distinct isomorphism types, cf.~Mendelsohn \& Hung~\cite{HuMe1972}, 1972. For $v=16$ there are exactly $1{,}054{,}163$ distinct isomorphism types, see Kaski, \"{O}sterg{\aa}rd \& Pottonen~\cite{KOP2006}, 2006.

If $\D=(X,\B,I)$ is a \mbox{$t$-$(v,k,\lambda)$} design with $t \geq
2$, and $x \in X$ arbitrary, then the \emph{derived} design with
respect to $x$ is \mbox{$\D_x=(X_x,\B_x, I_x)$}, where $X_x = X
\backslash \{x\}$, \mbox{$\B_x=\{B \in \B: (x,B)\in I\}$} and $I_x=
I \!\!\mid _{X_x \times \; \B_x}$. In this case, $\D$ is also called
an \emph{extension} of $\D_x$. Obviously, $\D_x$ is a
\mbox{$(t-1)$-$(v-1,k-1,\lambda)$} design.

For a group \mbox{$G \leq \Aut(\mathcal{D})$} of automorphisms of $\mathcal{D}$, let
$G_B$ denote the setwise stabilizer of a block $B \in \mathcal{B}$. All other notation is standard.

A detailed account on Steiner systems can be found, e.g., in~\cite{BJL1999}
and~\cite{ColMath2007}. Comprehensive survey articles on Steiner quadruple systems include~\cite{HartPh1992} and~\cite{LindRo1978}.


\section{Proof of the Main Theorem}\label{proof}

Let $\mathcal{D}$ be a non-trivial Steiner quadruple system \mbox{SQS$(v)$} of order $v$, and $N \leq G \leq \Aut(N)$ with a projective linear simple group $N=L_n(q)$, $n \geq 2$, $q$ a prime power. Here, $(n,q) \neq (2,2)$, $(2,3)$.
We consider the natural action of $G$ on the projective space $PG(n-1,q)$, $v=\frac{q^n-1}{q-1}$ (note that $L_n(q)$ has two doubly transitive permutation representations of the given degree if \mbox{$n>2$}). We remark that besides the fact that $N=L_n(q)$ is not simple for $(n,q) = (2,2)$,  $(2,3)$, it is obviously not possible to obtain in these cases a non-trivial \mbox{SQS$(v)$} by definition.

\begin{lem}\label{lem1}
Let $N=L_2(q)$, $v=q+1$. If $N \leq G \leq \aut(N)$ is \mbox{$3$-homogeneous}, then the cases as in Examples~\emph{\ref{Ex1}} and~\emph{\ref{Ex2}} hold in the Main Theorem.
\end{lem}

\begin{proof}
If $G$ is \mbox{$3$-homogeneous}, then in particular $G$ is flag-transitive. Hence, we may argue as in the corresponding case in~\cite{Hu2001} to obtain the known classes of Steiner quadruple systems. We note that in doing so, we only need to rely on~\cite{Deletal86}, not on the classification of the finite simple groups.
\end{proof}

\begin{lem}\label{lem2}
Let $N=L_2(q)$, $v=q+1$. If $N \leq G \leq \aut(N)$ is not \mbox{$3$-homogeneous}, then the case described in Example~\emph{\ref{Ex3}} holds in the Main Theorem.
\end{lem}

\begin{proof}
If we assume that $G$ is not \mbox{$3$-homogeneous}, which is the case if and
only if \mbox{$q \equiv 1$ (mod $4$)}, then $G$ has more than one orbit on the blocks and hence there cannot exist any flag-transitive \mbox{SQS$(v)$}. However, we will show that there exists a class of Steiner quadruple systems on which $G$ operates point \mbox{$2$-transitively}.
Since $N$ is a \mbox{$2$-transitive} permutation group, we may restrict ourselves to the case $N=G$.
As $PGL_2(q)$ is \mbox{$3$-homogeneous}, the unique orbit
under $PGL_2(q)$ on the \mbox{$3$-subsets} splits under $G$ in
exactly two orbits of equal length. By the definition of Steiner
quadruple systems, it follows that $G$ has exactly two orbits on
the blocks. These have equal length as for any block $B \in \B$ in
each orbit, the representation $G_B \rightarrow$ \mbox{Sym$(B) \cong
S_4$} is faithful and thus
\[G_B \cong S_4.\] We remark that $G_B$ has then four Sylow \mbox{$3$-subgroups}. By
Proposition~\ref{Hanani} and the fact
that \mbox{$q \equiv 1$ (mod $4$)}, we have to distinguish the
following two cases:

\smallskip
\indent Case (a): $q=3^{2d}$, $d \geq 1$.
\smallskip

Since $3 \mid q$, each Sylow \mbox{$3$-subgroup} has exactly one fixed point.
Thus, we have at most one orbit of length $4$ under $G_B$. On the
other hand, the normalizer of a Sylow \mbox{$3$-subgroup} in the symmetric
group $S_3$ is $S_3$ itself, hence $S_3$ fixes the respective fixed
point. The stabilizer of that fixed point in $S_4$ has order at least $6$.
But, as it is \mbox{$3$-closed}, it cannot be
$S_4$ itself. Moreover, it cannot be the alternating group $A_4$
because the latter does not contain $S_3$. Thus, it can only have
order $6$. Therefore, there exists at least one orbit of length $4$.
Hence, we have in each of the two orbits on the blocks exactly one
orbit of length $4$ under $G_B$. This yields the circle geometries described in Example~\ref{Ex3},
where we choose $a \not \in ({\mathbb{F}_q^*)}^2$, since
in general $-1 \in ({\mathbb{F}_q^*)}^2 \Leftrightarrow q \equiv 1$ (mod
$4$). As $2^4 \mid (3^d-1)(3^d+1)(3^{2d}+1)=3^{4d}-1$ for all $d \geq 1$,
we conclude that for $q=3^{2d}$, $d \geq 1$, always $q^2 \equiv
1$ (mod $16$) holds. Thus, we have in $G$
\[\frac{(q+1)q(q-1)}{24}\]
many subgroups isomorphic to $S_4$ on two conjugacy classes of equal
length (cf.~\cite[p.\,285]{Dick1901}). As we have precisely
\[\frac{b}{2}=\frac{(q+1)q(q-1)}{2 \cdot 24}\]
circles on each orbit of blocks, we obtain no further \mbox{SQS$(v)$}.

\smallskip
\indent Case (b): $q \equiv 1$ (mod $12$).
\smallskip

Let us assume that we have an orbit of length $4$ under $G_B$. Then,
the stabilizer, say $U$, of a point in $G_B$ is isomorphic to $S_3$. In
$U$, we have a normal subgroup of order $3$, which has exactly two fixed
points as in particular $3 \mid q-1$. But, as these are left fixed by an
involution in $U$, clearly $U$ has two fixed points. On the other hand,
the stabilizer on two points in $L_2(q)$ is cyclic, which leads to a
contradiction as $S_3$ is non-Abelian. Hence, there cannot exist any \mbox{SQS$(v)$}
in this case.
\end{proof}

\begin{lem}\label{lem3}
Let $N=L_n(q)$, $n \geq 3$, $v=\frac{q^n-1}{q-1}$. Then $N \leq G \leq \aut(N)$ cannot act as a group of automorphisms on any \mbox{SQS$(v)$}.
\end{lem}

\begin{proof}
Here $\Aut(N)=P \mathit{\Gamma} L_n(q) \rtimes
\text{\footnotesize{$\langle$}} \iota_\beta
\text{\footnotesize{$\rangle$}}$, where $\iota_\beta$ denotes the
graph automorphism induced by the inverse-transpose map
$\beta:GL_n(q) \rightarrow GL_n(q)$, $x
\mapsto {^t(x^{-1})}$. If $n=3$, then $v=q^2+q+1$ is always odd, in contrast to Proposition~\ref{Hanani}.
For $n>3$, we establish the claim via induction over $n$. Let us assume that there
is a counter-example with $n$ minimal. Without restriction, we can
choose three distinct points $x,y,z$ from a hyperplane $\h$ of
\mbox{$PG(n-1,q)$}. The translation group $T(\h)$ acts regularly on the points of
\mbox{$PG(n-1,q) \setminus \h$} and trivially on $\h$.
Hence, the unique block $B \in \B$ which is incident with the \mbox{$3$-subset}
$\{x,y,z\}$ must be contained completely in $\h$, since otherwise it would
contain all points of $PG(n-1,q) \setminus \h$, yielding a contradiction.
Thus, $\h$ induces a
\mbox{SQS$(\frac{q^{n-1}-1}{q-1})$} on which $G$ with socle $L_{n-1}(q)$
operates. Inductively, we obtain the minimal
counter-example for $n=3$, which we know is impossible. This verifies the claim.
\end{proof}

\noindent \textbf{Proof of Main Theorem}: The result is obtained by putting together Lemmas~\ref{lem1},~\ref{lem2}, and~\ref{lem3}.

\bibliographystyle{amsplain}
\bibliography{XbibPSL_SQS}
\end{document}